\documentclass[10pt,reqno]{amsart}

\usepackage{amsmath,amssymb,amsthm}

\usepackage{mathrsfs}

\usepackage[hidelinks]{hyperref}

\theoremstyle{plain}

\newtheorem{theorem}{Theorem}[section]

\newtheorem{proposition}[theorem]{Proposition}

\newtheorem{corollary}[theorem]{Corollary}

\theoremstyle{remark}

\newtheorem{remark}[theorem]{Remark}

\numberwithin{equation}{section}

\newcommand{\Ric}{\operatorname{Ric}}

\newcommand{\Hess}{\operatorname{Hess}}

\newcommand{\tr}{\operatorname{tr}}

\newcommand{\divg}{\operatorname{div}}

\newcommand{\Vol}{\operatorname{Vol}}

\newcommand{\Ein}{\operatorname{Ein}}

\newcommand{\Ad}{\operatorname{Ad}}

\newcommand{\Ker}{\operatorname{ker}}

\newcommand{\Img}{\operatorname{Im}}

\newcommand{\Span}{\operatorname{span}}

\begin{document}

\title[Harmonic critical metrics and Ricci flows]{Three-dimensional harmonic critical metrics for the Einstein-Hilbert functional and self-similar Ricci flows}

\author{Sergey E. Stepanov}

\address{Department of Mathematics, Financial University under the Government of the Russian Federation, 49 Leningradsky Prospect, 125993 Moscow, Russia}

\email{s.e.stepanov@mail.ru}

\author{Irina I. Tsyganok}

\address{Department of Mathematics, Financial University under the Government of the Russian Federation, 49 Leningradsky Prospect, 125993 Moscow, Russia}

\email{i.i.tsyganok@mail.ru}



\begin{abstract}

We address the three-dimensional harmonic analogue of Hilbert's variational problem: does criticality of the Einstein-Hilbert functional under volume-preserving Chen-Nagano variations force constant sectional curvature, and how is this rigidity related to Ricci flow and topology? Starting from the known characterization of compact harmonic critical metrics as gradient Ricci solitons, we use Ivey's classification to show that such metrics on closed three-manifolds are precisely the space-form metrics, equivalently the stationary metrics of the volume-normalized Ricci flow. We also describe residual gauge one-forms, establish transversality of the Chen-Nagano gauge to infinitesimal diffeomorphism directions, and give explicit self-similar flows. The normalized Ricci-flow velocity is an admissible Chen-Nagano variation for $n\ge 3$, and the smooth Poincar\'e theorem admits an equivalent formulation in terms of the existence of harmonic critical metrics.

\end{abstract}

\maketitle

\section{Introduction}\label{sec:intro}

The classical Hilbert variational problem characterizes Einstein metrics as critical points of the Einstein-Hilbert functional under fixed-volume variations. In its harmonic counterpart, the admissible variations are restricted by the Chen-Nagano condition. The general variational characterization of compact harmonic critical metrics is already available in dimensions $n \ge 3$~\cite{StepanovTsyganok2026b,Stepanov2027}. The problem of the present paper is therefore specifically three-dimensional: must such a critical metric have constant sectional curvature, and how is the resulting rigidity reflected in normalized Ricci flow and the topology of simply connected three-manifolds?

Let $(M^n, g)$ be a connected compact Riemannian manifold, $n \ge 3$, and set

\begin{equation}\label{eq:def_S_B}
\mathcal{S}_g = \int_M s_g \, dv_g, \qquad B_g(h) = \delta_g h + \frac{1}{2} d(\tr_g h), \qquad (\delta_g h)_j = -\nabla^i h_{ij}.
\end{equation}

Following Chen and Nagano~\cite{ChenNagano1984}, the condition $B_g(h) = 0$ is the harmonic symmetric-tensor condition. The admissible first-order volume-preserving directions form

\begin{equation}\label{eq:def_Vg}
\mathcal{V}_g = \left\{ h \in C^\infty(S^2 T^*M) : B_g(h) = 0, \quad \int_M \tr_g h \, dv_g = 0 \right\}.
\end{equation}

\begin{remark}[Historical background]\label{rem:history}

The differential condition $B_g(h) = 0$ links three distinct developments in differential geometry and mathematical physics. Its gauge-theoretic antecedent is de Donder's harmonic-coordinate condition in general relativity~\cite{deDonder1921}, whose linearization is closely related to the Bianchi gauge. Hamilton initiated the theory of Ricci flow and proved short-time existence using the Nash-Moser method~\cite{Hamilton1982}. DeTurck subsequently obtained a gauge-modified, strictly parabolic formulation~\cite{DeTurck1983}. In a distinct, intrinsic approach, Chen and Nagano~\cite{ChenNagano1984} studied harmonic metrics through the harmonicity of the identity map $\mathrm{id}\colon (M, g) \to (M, G)$. For a fixed background metric $g$, that condition is equivalently $B_g(G) = 0$; its infinitesimal version at $G = g$ is $B_g(h) = 0$. Chen and Nagano also showed that the space of harmonic symmetric two-tensors $\mathcal{H}_g=\Ker B_g$ is infinite-dimensional; see also~\cite{Chen2020}.

In the present variational problem this invariant condition restricts admissible metric variations, rather than serving as a gauge fixing for an evolution equation. Finally, the contracted second Bianchi identity gives $B_g(\Ric_g) = 0$, which motivates Proposition~\ref{prop:flow_admissible}. For further background on harmonic identity maps, see~\cite{Chen2020}.

\end{remark}

The general variational theory was developed in~\cite{Stepanov2026a,StepanovTsyganok2026a,StepanovTsyganok2026b,Stepanov2027}. We recall the precise characterization established in~\cite{StepanovTsyganok2026b,Stepanov2027}, which serves as the starting point for the present three-dimensional analysis.

\begin{theorem}[Harmonic variational characterization]\label{thm:char}

Let $(M^n, g)$ be a connected closed Riemannian manifold, $n \ge 3$. Then $g$ is critical for the Einstein-Hilbert functional with respect to every first-order volume-preserving Chen-Nagano variation if and only if $(M^n, g)$ is a gradient Ricci soliton. Equivalently, there exist a smooth function $f$ and a constant $\lambda$ such that $\Ric_g + \Hess_g f = \lambda g$. Necessarily, $\lambda = \bar{s}_g / n$.

\end{theorem}

\begin{proof}

We reproduce the variational argument of~\cite{StepanovTsyganok2026b,Stepanov2027} to make the present paper self-contained. Write $\mathcal{S}_g = \int_M s_g \, dv_g$ and let

\[
B_g(h) = \delta_g h + \frac{1}{2} d(\tr_g h), \qquad B_g^* \theta = \delta_g^* \theta + \frac{1}{2}(\delta_g \theta) g,
\]

where $(\delta_g h)_j = -\nabla^i h_{ij}$ and $\delta_g^* \theta = \frac{1}{2}\mathcal{L}_{\theta^\sharp} g$. Set $\mathcal{H}_g = \Ker B_g$ and $\mathcal{V}_g = \mathcal{H}_g \cap \{g\}^\perp$, where orthogonality is taken in $L^2(M^n,g)$. The principal symbol of $B_g^*$ is injective for $n \ge 3$; hence $B_g B_g^*$ is elliptic. On a closed manifold, elliptic regularity yields the $L^2$-orthogonal decomposition

\[
C^\infty(S^2 T^*M) = \Img B_g^* \oplus_{L^2} \Ker B_g.
\]

Since $B_g(g) = 0$, we also have $\mathcal{H}_g = \mathcal{V}_g \oplus \Span\{g\}$ orthogonally. The first variation of $\mathcal{S}$ is

\[
\mathcal{S}'_g(h) = -\int_M \langle \Ein_g, h \rangle \, dv_g, \qquad \Ein_g = \Ric_g - \frac{1}{2} s_g g.
\]

Consequently, $\mathcal{S}'_g(h) = 0$ for every $h \in \mathcal{V}_g$ if and only if there are a smooth one-form $\theta$ and a real constant $c$ for which

\[
\Ein_g = B_g^* \theta + c g.
\]

Taking the trace, using $\tr_g \Ein_g = -\frac{n-2}{2} s_g$ and $\tr_g(B_g^* \theta) = \frac{n-2}{2}\delta_g \theta$, and integrating over $M$, we obtain

\[
c = -\frac{n-2}{2n}\bar{s}_g, \qquad \bar{s}_g = \frac{1}{\Vol(M,g)}\int_M s_g \, dv_g.
\]

Here $\int_M \delta_g \theta \, dv_g = 0$ by Stokes' theorem. Since $n \ge 3$, the pointwise trace of the resulting Euler-Lagrange equation also gives $s_g + \delta_g \theta = \bar{s}_g$. Substitution into the tensor equation now yields

\[
\Ric_g = \delta_g^* \theta + \frac{\bar{s}_g}{n} g.
\]

Putting $X = -\theta^\sharp$ and $\lambda = \bar{s}_g / n$, we conclude that

\[
\Ric_g + \frac{1}{2}\mathcal{L}_X g = \lambda g.
\]

Thus $(M^n, g, X)$ is a compact Ricci soliton. By Perelman's gradient characterization of compact Ricci solitons~\cite{Perelman2002}, there exists a smooth function $f$ such that $\Ric_g + \Hess_g f = \lambda g$.

Conversely, suppose $\Ric_g + \Hess_g f = \lambda g$ on a closed manifold. Tracing and integrating gives $\lambda = \bar{s}_g / n$. Set $\theta = -df$. Then $\delta_g^* \theta = -\Hess_g f$ and $\delta_g \theta = \Delta_g f$, where $\Delta_g f = \divg_g(\nabla f)$. The soliton trace identity is $s_g + \Delta_g f = \bar{s}_g$, and therefore

\[
\Ein_g = B_g^* \theta - \frac{n-2}{2n}\bar{s}_g g.
\]

For each $h \in \mathcal{V}_g$, the adjoint relation and the volume constraint imply

\[
\int_M \langle B_g^* \theta, h \rangle \, dv_g = \int_M \langle \theta, B_g(h) \rangle \, dv_g = 0, \qquad \int_M \tr_g h \, dv_g = 0.
\]

Hence $\mathcal{S}'_g(h) = 0$ for every admissible $h$, proving harmonic criticality. Moreover, if $\theta$ is any one-form satisfying the Euler-Lagrange equation, then $\delta_g^*(\theta + df) = 0$; thus $\theta = -df + \eta$, where $\eta^\sharp$ is a Killing vector field.

\end{proof}

Theorem~\ref{thm:char} holds in every dimension $n \ge 3$, whereas the rigidity result below relies on Ivey's classification of compact three-dimensional Ricci solitons~\cite{Ivey1993}. We apply these results to the harmonic variational framework to describe the geometry, residual gauge fields, and explicit Ricci-flow evolution, including a transversality consequence for the Chen-Nagano gauge. The general first-variation argument has been reproduced above from~\cite{StepanovTsyganok2026b,Stepanov2027}. We additionally identify the normalized Ricci-flow velocity as an admissible harmonic variation and give a precise variational reformulation of the smooth Poincar\'e theorem.

\section{Three-dimensional variational rigidity}\label{sec:rigidity}

The three-dimensional case is distinguished from the higher-dimensional theory. Already in real dimension four, compact non-Einstein gradient shrinking K\"ahler-Ricci solitons exist (the Koiso-Cao example; see~\cite{Cao2006}); by Theorem~\ref{thm:char} they are harmonic critical but not critical in the classical fixed-volume Hilbert problem. In dimension three, by contrast, the two variational problems have precisely the same critical metrics.

\begin{theorem}\label{thm:rigidity}
Let $(M^3, g)$ be a connected closed Riemannian manifold. The following conditions are equivalent:
\begin{enumerate}
\item[(i)] $g$ is critical for $\mathcal{S}$ in every direction $h \in \mathcal{V}_g$;
\item[(ii)] $g$ admits a (necessarily gradient) compact Ricci-soliton structure;
\item[(iii)] $g$ is an Einstein metric;
\item[(iv)] $g$ has constant sectional curvature $\kappa \in \mathbb{R}$;
\item[(v)] $g$ is a stationary solution of the volume-normalized Ricci flow.
\end{enumerate}
Thus restricted harmonic criticality and ordinary fixed-volume criticality coincide for compact three-manifolds. Moreover, $\lambda = 2\kappa$ and $s_g = \bar{s}_g = 6\kappa$.
\end{theorem}

\begin{proof}
The equivalence $\mathrm{(i)} \Leftrightarrow \mathrm{(ii)}$ follows from the restricted first-variation theorem in~\cite{Stepanov2027}, summarized by Theorem~\ref{thm:char}. By Ivey's classification~\cite{Ivey1993}, every compact three-dimensional Ricci soliton has constant sectional curvature, proving $\mathrm{(ii)} \Rightarrow \mathrm{(iv)}$. A space-form metric satisfies $\Ric_g = 2\kappa g$, hence $\mathrm{(iv)} \Rightarrow \mathrm{(iii)} \Rightarrow \mathrm{(ii)}$ with a constant soliton potential. Finally, on a three-manifold the Riemann curvature tensor is determined by $\Ric_g$, so $\mathrm{(iii)} \Rightarrow \mathrm{(iv)}$. The classical first-variation principle identifies Einstein metrics as the unrestricted fixed-volume critical points of $\mathcal{S}$~\cite{Besse1987}. For a constant-curvature metric, $-2\Ric_g + \frac{2}{3}\bar{s}_g g = 0$; conversely, the stationarity equation implies $\Ric_g = \frac{1}{3}\bar{s}_g g$, so $\mathrm{(iii)} \Leftrightarrow \mathrm{(v)}$. The numerical identities follow by taking traces.
\end{proof}

\begin{remark}\label{rem:thurston}
The shrinking condition $\lambda > 0$ is equivalent to $\kappa > 0$. Harmonic critical metrics on closed three-manifolds also include flat ($\lambda = 0$) and hyperbolic ($\lambda < 0$) space forms. No other Thurston model occurs as the constant-curvature geometry of such a metric. This statement concerns compact critical metrics and does not classify noncompact homogeneous Ricci solitons. By contrast, non-Einstein compact K\"ahler shrinking solitons exist already in real dimension four; see~\cite{Cao2006}.
\end{remark}

\section{Residual harmonic gauge}\label{sec:gauge}

The formal adjoint of $B_g$ is
\begin{equation}\label{eq:Bg_adjoint}
B_g^* \theta = \delta_g^* \theta + \frac{1}{2}(\delta_g \theta) g, \qquad \delta_g \theta = -\divg_g \theta^\sharp.
\end{equation}
The restricted Euler-Lagrange equation derived in~\cite{Stepanov2027} takes the form
\begin{equation}\label{eq:restricted_EL}
\Ein_g = B_g^* \theta - \frac{n-2}{2n}\bar{s}_g g, \qquad \Ein_g = \Ric_g - \frac{1}{2}s_g g.
\end{equation}

\begin{theorem}\label{thm:gauge}
For a compact three-dimensional harmonic critical metric, every one-form $\theta$ satisfying~\eqref{eq:restricted_EL} is dual to a Killing vector field, and conversely every Killing one-form satisfies~\eqref{eq:restricted_EL}. More precisely: if $\kappa < 0$, then $\theta = 0$; if $\kappa = 0$, the admissible one-forms form the $b_1(M)$-dimensional space of parallel one-forms; and if $\kappa > 0$, they form the finite-dimensional space of Killing one-forms of the spherical space form $S^3/\Gamma$, isomorphic to the centralizer of $\Gamma$ in $\mathfrak{so}(4)$ (of dimension at most $6$, with equality precisely for $S^3$ and $\mathbb{R}P^3$).
\end{theorem}

\begin{proof}
By Theorem~\ref{thm:rigidity}, $\Ric_g = 2\kappa g$ and $s_g = \bar{s}_g = 6\kappa$. Thus~\eqref{eq:restricted_EL} reduces to $B_g^* \theta = 0$. Taking the trace of~\eqref{eq:Bg_adjoint} in dimension three yields $\tr_g(B_g^* \theta) = \frac{1}{2}\delta_g \theta$, so $\delta_g \theta = 0$ and consequently $\delta_g^* \theta = 0$, which is the Killing equation for $\theta^\sharp$. Conversely, a Killing field has zero divergence and satisfies $B_g^* \theta = 0$. For a Killing field $K$ on a closed manifold, the Bochner identity implies

\[\int_M |\nabla K|^2 \, dv_g = \int_M \Ric_g(K, K) \, dv_g.\]

If $\kappa < 0$, the right-hand side is nonpositive, hence $K = 0$. If $\kappa = 0$, then $\nabla K = 0$, so $K^\flat$ is harmonic and spans a space of dimension $b_1(M) \in \{0, 1, 2, 3\}$ (in particular, $\theta = 0$ on the Hantzsche-Wendt manifold, where $b_1(M) = 0$). By Hodge theory and the Bochner formula, every harmonic one-form on a compact flat manifold is parallel. For $\kappa > 0$, the Killing algebra of $S^3/\Gamma$ is the centralizer of $\Gamma$ in $\mathfrak{so}(4)$.

To spell out the final group-theoretic step, the centralizer is the fixed-point subalgebra of the adjoint action of $\Gamma$ on $\mathfrak{so}(4)$. Since $\dim \mathfrak{so}(4) = 6$, a six-dimensional centralizer means $\Ad(\gamma) = \mathrm{id}$ for every $\gamma \in \Gamma$. Hence $\Gamma \subset \ker \Ad = Z(\mathrm{SO}(4)) = \{\pm I\}$. The corresponding spherical space forms are $S^3$ and $\mathbb{R}P^3$.
\end{proof}

The relation between infinitesimal harmonic transformations and the Yano
operator was investigated in~\cite{StepanovShandra2003,StepanovTsyganokMikes2013}.
The preceding theorem describes the one-forms occurring in the restricted
Euler-Lagrange equation. A related, but distinct, assertion concerns
pure-gauge deformations of the metric.

\begin{corollary}[Transversality of the Chen-Nagano gauge]
\label{cor:gauge_transversality}
Let $(M^3,g)$ be a connected closed Riemannian three-manifold with a
harmonic critical metric. For every smooth vector field $X$ on $M^3$,
the condition $B_g(\mathcal L_Xg)=0$ is equivalent to $X$ being a
Killing vector field (that is, $\mathcal L_Xg=0$).
Consequently,
\[
 \mathcal V_g\cap\{\mathcal L_Xg:X\in\Gamma(TM)\}=\{0\}.
\]
In particular, the Chen-Nagano gauge contains no nonzero infinitesimal
pure-gauge metric deformations at a three-dimensional harmonic critical metric.
\end{corollary}

\begin{proof}
By Theorem~\ref{thm:rigidity}, the metric $g$ has constant sectional
curvature $\kappa$, so $\Ric_g=2\kappa g$. Put $\theta=X^\flat$.
By the Yano-operator identity for infinitesimal harmonic
transformations established in~\cite{StepanovShandra2003}
(see also~\cite{StepanovTsyganokMikes2013}),
\[
 B_g(\mathcal L_Xg)=\square_g\theta,
 \qquad \square_g=\Delta_H-4\kappa I,
\]
where $\Delta_H=d\delta_g+\delta_gd$ is the nonnegative Hodge-de Rham
Laplacian. Here $\mathcal L_Xg=2\delta_g^*\theta$, consistently with
our normalization of $\delta_g^*$.

Suppose that $B_g(\mathcal L_Xg)=0$. The Hodge decomposition gives
\[
 \theta=du+\eta,\qquad \delta_g\eta=0,
 \qquad\int_Mu\,dv_g=0.
\]
Since $\square_g=\Delta_H-4\kappa I$ preserves exact and co-closed
one-forms, the orthogonality of these two components yields
$\square_g(du)=0$ and $\square_g\eta=0$.
The exact component therefore satisfies
\[
 d(\Delta_gu+4\kappa u)=0,
\]
where $\Delta_g=\operatorname{div}_g\nabla$ is the scalar Laplacian
with nonpositive spectrum. By the mean-zero normalization,
\[
 (\Delta_g+4\kappa)u=0.
\]
If $\kappa\leq 0$, integration by parts implies $u=0$.
If $\kappa>0$, lifting $u$ to the round universal covering sphere
of curvature $\kappa$ gives a solution of the same equation.
The eigenvalues of $-\Delta_g$ on that sphere are
$\kappa\ell(\ell+2)$, $\ell=0,1,2,\ldots$; none equals $4\kappa$.
Thus $u=0$ in this case as well.

We conclude that $\theta=\eta$, $\delta_g\theta=0$, and
$B_g(\mathcal L_Xg)=0$. By the adjoint relation and
$B_g^*\theta=\frac12\mathcal L_Xg$, it follows that
\[
 0=\int_M\langle B_g(\mathcal L_Xg),\theta\rangle\,dv_g
 =\frac12\int_M|\mathcal L_Xg|^2\,dv_g.
\]
Consequently $\mathcal L_Xg=0$, and the converse is immediate.
Finally,
\[
 \int_M\operatorname{tr}_g(\mathcal L_Xg)\,dv_g
 =2\int_M\operatorname{div}_gX\,dv_g=0.
\]
Thus every pure-gauge variation on a closed manifold satisfies the
first-order volume constraint, and the asserted intersection is trivial.
\end{proof}

\section{Self-similar Ricci flow and topology}\label{sec:flow_topology}

\begin{proposition}[Ricci-flow directions are admissible Chen-Nagano variations]\label{prop:flow_admissible}
Let $(M^n, g(t))$, $n \ge 3$, be a smooth solution of the volume-normalized Ricci flow on a connected closed manifold. At every time $t$ for which the solution exists, its instantaneous velocity

\[h_g := \partial_t g = -2\Ric_g + \frac{2}{n}\bar{s}_g g\]

belongs to the admissible space $\mathcal{V}_g$ defined in~\eqref{eq:def_Vg}. Equivalently,

\[B_g(h_g) = 0, \qquad \int_M \tr_g h_g \, dv_g = 0.\]

\end{proposition}

\begin{proof}
With our divergence convention $(\delta_g h)_j = -\nabla^i h_{ij}$, the contracted Bianchi identity gives

\[B_g(\Ric_g) = \delta_g \Ric_g + \frac{1}{2} d s_g = 0.\]

Since $B_g(g) = 0$ and $\bar{s}_g$ is spatially constant, the definition of $h_g$ implies $B_g(h_g) = 0$. Furthermore,

\[\int_M \tr_g h_g \, dv_g = -2\int_M s_g \, dv_g + 2\bar{s}_g \Vol(M,g) = 0.\]

Thus $h_g \in \mathcal{V}_g$ at every time under consideration.
\end{proof}

The admissible space $\mathcal{V}_{g(t)}$ varies with the evolving metric. Consequently, Proposition~\ref{prop:flow_admissible} concerns the instantaneous tangent direction of the normalized flow; it does not identify that flow with a gradient flow of the Einstein-Hilbert functional. In dimension three, Theorem~\ref{thm:rigidity} separately identifies its stationary metrics with the harmonic critical metrics.

\begin{corollary}\label{cor:explicit_flow}
Let $(M^3, g)$ satisfy Theorem~\ref{thm:rigidity} and have sectional curvature $\kappa$. The unique unnormalized Ricci flow starting at $g$ is
\begin{equation}\label{eq:unnorm_flow}
g(t) = (1 - 4\kappa t)g, \qquad s_{g(t)} = \frac{6\kappa}{1 - 4\kappa t}
\end{equation}
for the times when $1 - 4\kappa t > 0$. If $\kappa > 0$, extinction occurs at $T = (4\kappa)^{-1}$; if $\kappa = 0$, the solution is stationary; and if $\kappa < 0$, it expands for all $t \ge 0$. Under the volume-normalized flow
\begin{equation}\label{eq:norm_flow}
\partial_t g = -2\Ric_g + \frac{2}{3}\bar{s}_g g,
\end{equation}
$g(t) = g$ for all $t \ge 0$.
\end{corollary}

\begin{proof}
Since $\Ric_g = 2\kappa g$ and the Ricci tensor as a $(0,2)$-tensor is invariant under a spatially constant rescaling of the metric, the right-hand side of the Ricci flow is $-4\kappa g$, which gives~\eqref{eq:unnorm_flow}. The formula for scalar curvature follows from scaling. Since $\bar{s}_g = 6\kappa$, the two terms in~\eqref{eq:norm_flow} cancel. The asserted uniqueness is the standard uniqueness of smooth Ricci flow on a compact manifold.
\end{proof}

\begin{corollary}\label{cor:finite_pi1}
Let $(M^3,g)$ be a connected closed harmonic critical three-manifold with finite fundamental group. Then $g$ has positive sectional curvature ($\kappa > 0$) and $(M^3, g)$ is isometric, up to constant rescaling, to a spherical space form $S^3/\Gamma$ with $\Gamma$ finite. In particular, if $M^3$ is simply connected, then $(M^3, g)$ is isometric to the round three-sphere $S^3$ up to scaling. On a fixed compact hyperbolic three-manifold, any two negatively curved harmonic critical metrics are isometric up to constant rescaling.
\end{corollary}

\begin{proof}
Every complete simply connected manifold of nonpositive sectional curvature is noncompact by Cartan-Hadamard; therefore a compact space form with finite fundamental group must have $\kappa > 0$. Up to rescaling $\kappa$ to $1$, its universal cover is isometric to the unit round sphere $S^3$, so $(M^3, g)$ is isometric to a quotient $S^3/\Gamma$ by a finite group of isometries acting freely. In the simply connected case $\Gamma = \{1\}$, yielding the round metric on $S^3$ up to scaling. For a compact hyperbolic three-manifold, normalizing both curvatures to $-1$ and applying Mostow rigidity~\cite{Mostow1968} yields the last assertion.
\end{proof}

\begin{corollary}[A variational formulation of the smooth Poincar\'e theorem]\label{cor:poincare}
Let $M^3$ be a connected closed simply connected smooth manifold. The following conditions are equivalent:
\begin{enumerate}
\item[(i)] $M^3$ is diffeomorphic to $S^3$;
\item[(ii)] $M^3$ admits a Riemannian metric $g$ that is critical for the Einstein-Hilbert functional with respect to all $h \in \mathcal{V}_g$;
\item[(iii)] $M^3$ admits a Riemannian metric of constant positive sectional curvature.
\end{enumerate}
Whenever $\mathrm{(ii)}$ holds, every harmonic critical metric on $M^3$ is isometric, up to a constant rescaling, to the round metric on $S^3$. Consequently, the universal assertion that every connected closed simply connected smooth three-manifold carries a harmonic critical metric is logically equivalent to the smooth three-dimensional Poincar\'e theorem.
\end{corollary}

\begin{proof}
The implication $\mathrm{(i)} \Rightarrow \mathrm{(iii)}$ follows by pulling back the round metric on $S^3$ through a diffeomorphism. A metric in $\mathrm{(iii)}$ is Einstein and therefore harmonic critical by Theorem~\ref{thm:rigidity}, which yields $\mathrm{(iii)} \Rightarrow \mathrm{(ii)}$. Finally, Corollary~\ref{cor:finite_pi1} shows that $\mathrm{(ii)}$ forces the metric to be round, up to scaling, proving $\mathrm{(ii)} \Rightarrow \mathrm{(i)}$, as well as the asserted rigidity of every harmonic critical metric. The universal equivalence follows immediately.
\end{proof}

\begin{remark}[Relation to Perelman's theory and logical scope]\label{rem:scope}
Corollary~\ref{cor:poincare} does not provide an independent proof of the smooth Poincar\'e theorem. The implication from harmonic criticality to spherical geometry follows from Corollary~\ref{cor:finite_pi1} and Ivey's classification~\cite{Ivey1993}; no general existence theorem for harmonic critical metrics on arbitrary closed simply connected three-manifolds is proved here.

The comparison with Perelman's variational framework is nevertheless instructive. Perelman introduced entropy functionals depending on the metric $g$ and an auxiliary potential $f$ (and, in the case of his $\mathcal{W}$-functional, a scale parameter $\tau$). His $\mathcal{F}$- and $\mathcal{W}$-functionals are naturally associated with gradient steady and gradient shrinking Ricci solitons, respectively, and yield monotonicity formulas for appropriately coupled Ricci flows~\cite{Perelman2002}. In contrast, the Chen-Nagano approach retains the classical Einstein-Hilbert functional $\mathcal{S}_g$ and restricts the admissible infinitesimal metric variations to $\mathcal{V}_g$. By Theorem~\ref{thm:char}, its compact critical metrics are precisely the compact gradient Ricci soliton metrics, including the Einstein metrics of positive, zero, or negative scalar curvature. Thus these approaches meet in Ricci-soliton geometry, but their functionals and critical sets are not identical in general.

In dimension three, Ivey's classification reduces all compact solitons to metrics of constant sectional curvature. The Hamilton-Perelman theory applies to arbitrary initial Riemannian metrics and establishes the smooth Poincar\'e theorem using Ricci flow with surgery and finite-time extinction~\cite{Perelman2003a,Perelman2003b}. Once this topological theorem is known, pulling back the round metric on $S^3$ gives the harmonic critical metric required by condition~\textup{(ii)} of Corollary~\ref{cor:poincare}. Our corollary is therefore a conditional geometric rigidity result and a logical variational reformulation, not an alternative existence proof.
\end{remark}

\end{document}